\documentclass{amsart}
\usepackage{color}
\usepackage{diagrams}

\usepackage{amssymb}
\usepackage{amsmath}
\usepackage{amsthm}
\usepackage{amsfonts}
\usepackage{amstext}
\usepackage{amsopn}
\usepackage{amscd}

\newtheorem{theorem}{Theorem}[section]

\newtheorem{corollary}[theorem]{Corollary}

\newtheorem{lemma}[theorem]{Lemma}

\newtheorem{proposition}[theorem]{Proposition}
\newtheorem{question}[theorem]{Question}

\newtheorem*{cor_closed}{Corollary \ref{closedn}}

\newtheorem*{thm_sym}{Theorem \ref{symplectic}}
\newtheorem*{thm_2complex}{Theorem \ref{2complex}}
\newtheorem*{prop_def}{Proposition \ref{defG}}

\theoremstyle{definition}
\newtheorem{definition}[theorem]{Definition}
\newtheorem{example}[theorem]{Example}

\newtheorem{remark}[theorem]{Remark}

\newcommand{\db}{\bar{\delta}}
\newcommand{\R}{\operatorname{r}}
\newcommand{\Z}{\mathbb{Z}}
\newcommand{\K}{\mathbb{K}}
\newcommand{\rk}{\operatorname{rank}}
\newcommand{\xp}{X^\prime}
\newcommand{\onto}{\twoheadrightarrow}
\newcommand{\Ext}{\operatorname{Ext}}

\newcommand{\Def}{\operatorname{def}}

\begin{document}

\title[Monotonicity of the Degrees]{Monotonicity of Degrees of Generalized Alexander Polynomials of Groups and 3-Manifolds}

\author[Shelly L. Harvey]{Shelly L. Harvey$^\dag$}
\address{Department of Mathematics, MIT \\
Cambridge, MA, 02319-4307}
\email{sharvey@math.mit.edu}

\thanks{$^{\dag}$The author was partially supported by an NSF Postdoctoral Fellowship}
\date{}

\begin{abstract} We investigate the behavior of the higher-order degrees, $\db_n$, of a finitely presented group $G$. 
These $\db_n$ are functions from $H^1(G;\Z)$ to $\Z$ whose values are the degrees certain higher-order
Alexander polynomials.  We show that if $\Def(G)\geq 1$ or $G$ is the fundamental group of a compact, 
orientable 3-manifold then $\db_n$ is a monotonically increasing function of $n$ for $n\geq1$.  This is false
for general groups.
As a consequence, we show that 
if a 4-manifold of the form $X \times S^1$ admits a symplectic structure then X ``looks algebraically like'' a 
3-manifold that fibers over $S^1$, supporting a positive answer to a question of Taubes.  This generalizes a theorem of
S. Vidussi \cite{Vi} and is an improvement on the results in \cite{Ha1}.  We also find new conditions on a 3-manifold $X$ which will guarantee that the Thurston norm of $f^{\ast}(\psi)$, for $\psi \in H^1(X;\Z)$ and $f:Y \rightarrow X$ a surjective map on $\pi_1$, will be at least as large the Thurston norm of $\psi$.  When $X$ and $Y$ are knot complements, this gives a partial answer to a question of J. Simon.

More generally, we define $\Gamma$-degrees, $\db_\Gamma$, corresponding to a surjective map $G\twoheadrightarrow \Gamma$ for which $\Gamma$ is poly-torsion-free-abelian.  Under certain conditions, we show they satisfy a monotonicity condition if one varies the group.  As a result, we show that these generalized degrees  give obstructions to the
deficiency of a group being positive and obstructions to a finitely presented group being the fundamental group
of a compact, orientable 3-manifold.

\end{abstract}

\maketitle


In \cite{Ha1}, we defined some new invariants  $\db_n$ for a finite CW-complex $X$.  These invariants depended only on the fundamental group of $X$ and measured the ``size''
of the successive quotients of the rational derived series of $\pi_1(X)$.  Given $X$ and a
cohomology class $\psi \in H^1(X)$, $\db_n(\psi)$ was defined to be the degree of a ``higher-order Alexander polynomial.''
 Although defined algebraically, these degrees have many topological applications in the case that $X$ is
  a 3-manifold.  
  In this case, we showed that the $\db_n$ give new estimates for the Thurston norm of a 3-manifold generalizing a theorem of C. McMullen
  \cite{McM}.  
  Recall that the \textbf{Thurston norm} of a class $\psi \in H^1(X;\Z)$, $||\psi||_T$, is defined to be the minimum negative euler characteristic
of all (possibly disconnected) surfaces $F$ whose homology class $[F] \in H_2(X,\partial X;\Z)$ is Poincare dual to $\psi$ and such that each component of $F$ is
non-positively curved.
  The $\db_n$  also give new algebraic obstructions to a 3-manifold fibering over $S^1$, to a 4-manifold of
   the form $X \times S^1$ admitting a symplectic structure, and to a 3-manifold being
   Seifert fibered.   They were also shown to have applications to minimal ropelength and genera of knots and links in
   $S^3$.  Related work has been done by T. Cochran, K. Orr, P. Teichner, for knots and knot concordance in \cite{Co}
   and \cite{COT}.  Recently, V. Turaev  \cite{turaev} has generalized some of the results in \cite{Ha1}.

Since $\db_n$ only depends on the fundamental group, we can consider $\db_n$ as an invariant of a general group $G$, $\db_n: H^1(G;\Z) \rightarrow \Z$.  In this paper, we continue to investigate the special behavior of the $\db_n$ when $G$ is the fundamental group of a 3-manifold (with empty or toroidal boundary) or a group with deficiency at least 1. The results give new algebraic information about the topology of a symplectic
4-manifolds of the form $X \times S^1$.  They give obstructions to a finitely presented group having positive
deficiency or being the fundamental group of a compact, orientable 3-manifold (with or without boundary).  They also give new information about the behavior of the Thurston norm under a map between 3-manifolds which is surjective on $\pi_1$.  We state some of our main theorems and their applications below. 

In \cite{Ha1}, we constructed examples of 3-manifolds for which $\db_n$ was a strictly increasing function of $n$
for $n \geq 0$.  Moreover, it was conjectured that the $\db_n$ were always a monotonically increasing function of $n$
for $n \geq 1$.  We show that this conjecture is true.  By $\db_n \leq \db_{n+1}$ (respectively $\db_n=0$) we mean that
$\db_n(\psi) \leq \db_{n+1}(\psi)$ (respectively $\db_n(\psi)=0$) for all $\psi \in H^1(X)$.
\begin{cor_closed}Let $X$ be a  closed, orientable, connected 3-manifold.    If $\beta_1(X)\geq 2$  then
\[\db_0 \leq \db_1 \leq \cdots \leq \db_n \leq \cdots.
\]If $\beta_1(X)=1$ and $\psi$ is a generator
of $H^1(X)$ then $\db_0(\psi) -2 \leq \db_1(\psi) \leq \cdots \leq \db_n(\psi) \leq \cdots$.
\end{cor_closed}

As a consequence of Corollary~\ref{closedn}, we show (in Theorem \ref{symplectic}) that if 4-manifold of the form $X \times S^1$ admits a symplectic
structure then $X$ ``looks algebraically like'' a 3-manifold which fibers over $S^1$, thus further supporting a conjecture
 of Taubes.
The proof of Theorem~\ref{symplectic} uses a theorem of Vidussi in \cite{Vi} who proves this theorem in the case $n=0$.

\begin{thm_sym}Let $X$ be an closed, irreducible
3-manifold such that $X \times S^1$ admits a symplectic structure.
If $\beta_1(X) \geq 2$ there exists a $\psi \in
H^1(X;\mathbb{Z})$ such that \[\db_0(\psi)=\db_1(\psi)= \cdots
\db_n(\psi) = \cdots =\|\psi\|_T.
\] If $\beta _{1}( X) =1$ then for
any generator $\psi$ of $H^1(X;\mathbb{Z})$,
\[
\db_0(\psi)-2=\db_1(\psi)=
\cdots
 \db_n(\psi) = \cdots =\|\psi\|_T.
\]\end{thm_sym}

More generally, we define $\db_\Gamma(\psi)$ for any group $G$ and
any ``admissible pair'' $(\phi_{\Gamma}: G \onto \Gamma, \psi:G \onto \Z)$ of $G$.
When $G$ is a finitely presented group with $\Def(G)\geq 1$, we show that the $\db_n$ satisfy a monotonicity condition.
We also prove a similar theorem when $G$ is the fundamental group of a closed, orientable 3-manifold (see Theorem~\ref{closed}).
\begin{thm_2complex}Let $G$ be a finitely presented group with $\Def(G)\geq 1$ and
$(\phi_{_\Lambda},\phi_{_\Gamma},
\psi)$ be an admissible triple for $G$.  If $(\phi_{_\Lambda},\phi_{_\Gamma},\psi)$ is not initial then
\begin{equation}\db_{\Lambda}(\psi) \geq \db_{\Gamma}(\psi)\end{equation} otherwise
\begin{equation} \db_{\Lambda}(\psi) \geq \db_{\Gamma}(\psi)-1.
\end{equation}
\end{thm_2complex}

%

As a consequence of the monotonicity theorems, we see that the
$\db_\Gamma$ give obstructions to the deficiency of a group being positive or being the fundamental group of a compact,
orientable 3-manifold.  These obstructions are non-trivial even when the groups $\Gamma$ and $\Lambda$ are abelian. For example, we can easily recover the (known) result that $\Z^m$ cannot be the fundamental group of a compact 3-manifold when $m \geq 4$ (see below or for more details see Example~\ref{Zm}).
\begin{prop_def}Let $G$ be a finitely presented group and $(\phi_{_\Lambda},\phi_{_\Gamma},\psi)$ be an admissible triple
for $G$.
\begin{enumerate}\item Suppose $(\phi_{_\Lambda},\phi_{_\Gamma},\psi)$ is not initial.
If $\db_\Lambda(\psi) < \db_\Gamma(\psi)$ then $\Def(G) \leq 0$ and $G$ cannot be the fundamental group of
a compact, orientable 3-manifold (with or without boundary).
\item Suppose $(\phi_{_\Lambda},\phi_{_\Gamma},\psi)$ is initial.
If $\db_\Lambda(\psi) < \db_\Gamma(\psi) - 1$ then $\Def(G) \leq 0$ and $G$ cannot be the fundamental group
of a compact, orientable 3-manifold with at least one boundary component which is not a 2-sphere.  In addition,
if $\db_\Lambda(\psi) < \db_\Gamma(\psi) - 2$ then $G$ cannot be the fundamental group of a compact, orientable 3-manifold
(with or without boundary). \end{enumerate}
\end{prop_def}
Let us consider the simplest case when $\Lambda$ is the abelianization (modulo torsion) of $G$ and $\Gamma=\Z$.  
In this case, $\db_\Z(\psi)$ is equal to the rank of $H_1((X_G)_\psi;\Z)$ as an abelian group where $(X_G)_\psi$ is the infinite cyclic cover of $X_G$, a finite CW-complex with $\pi_1(X_G)=G$, corresponding to $\psi$ (as long as this number is finite).
Moreover,  $\db_{\Lambda}(\psi)$ is equal the Alexander norm of $\psi$ which depends only on $\psi$ and the multivariable Alexander polynomial of $G$.  For example, then the Alexander polynomial of $\Z^m$ is 1 so $\db_{\Z^m}(\psi)=0$ for any $\psi$.  Moreover, the first homology of any infinite cyclic cover of the $m$-torus is $\Z^{m-1}$ so $\db_\Z(\psi)=m-1$.  Thus, as mentioned above, we see that $\Z^m$ cannot be the fundamental group of a compact 3-manifold.

Recall that the $i^{th}$-order degree of a group $\db_i(\psi)$ is a specific example of the degree $\db_\Gamma(\psi)$. We give examples of finite 2-complexes $X_{n,g}$ with $\beta_1(X_{n,g})=1$ for $n,g \geq 1$ such that the $i^{th}$-order
degrees for $0\leq i \leq n-1$ of $X_{n,g}$ are ``large'' but the $n^{th}$-order degree is 0.  Thus the fundamental group  of these spaces cannot have positive deficiency nor can they be the fundamental group of a compact, orientable 3-manifold (see Proposition~\ref{remark2} and Example~\ref{ex}).
 
 Theorem~\ref{closed} also has applications to the study of the behavior of the genus of a knot under a surjective map on $\pi_1$.  
The following question was asked by J. Simon (see R. Kirby's Problem List \cite[Question 1.12(b)]{Ki}). 

\theoremstyle{definition}
\newtheorem*{jsimon}{Question 1.12(b) of \cite{Ki}}
\theoremstyle{plain}
\begin{jsimon}[J. Simon] If $J$ and $K$ are knots in $S^3$ and $f: S^3\setminus L \rightarrow S^3 \setminus K$ is surjective
on $\pi_1$, is $g(L) \geq g(K)$?
\end{jsimon}
The answer to the above question is known to be ``yes'' when $\delta_0(K)=2 g(K)$.  We strengthen this result to the case when $\delta_{n}(K)=2 g(K)-1$.  

\newtheorem*{knotincreasing}{Corollary~\ref{knot_greater}}
\begin{knotincreasing}Suppose $J$ and $K$ are knots in $S^3$ such that there 
exists a surjective homomorphism $\rho: \pi_1(S^3\setminus L) \twoheadrightarrow \pi_1(S^3\setminus K)$.  If $\db_0(K)=2 g(K)$ or $\db_n(K)=2 g(K) -1$ for
some $n \geq 1$ then $g(L)\geq g(K)$. 
\end{knotincreasing}

We also prove this is the case if we replace the genus of a knot by the Thurston norm.  The  following corollary is a generalization of the result due to Gabai \cite{Ga} that a degree one map $f:X\rightarrow Y$ between three manifolds gives the inequality $||f^{\ast}(\psi)||_T \geq ||\psi||_T$ for all $\psi \in H^1(Y;\Z)$. 
For simplicity, we state only the case when $\beta_1(Y)\geq 2$.

\newtheorem*{th_gr}{Corollary~\ref{thurston_greater}}
\begin{th_gr} Suppose there exists an epimorphism $\rho: \pi_1(X) \onto \pi_1(Y)$, where $X$ and $Y$ are compact,
orientable $3$-manifolds, with toroidal or empty boundaries, such that $\beta_1(X)=\beta_1(Y)\geq 2$ and $\R_0(X)=0$.
Let $\psi \in H^1(\pi_1(Y);\Z)$. If $\db_n(\psi)=||\psi||_T$ for some $n \geq 0$
then $$||\rho^{\ast}(\psi)||_T \geq ||\psi||_T.$$
\end{th_gr}

\section{Definitions}\label{definitions}

We will define the higher-order degrees $\db_\Gamma$ and ranks $\R_{_\Gamma}$ of a group $G$ and surjective homomorphism $\phi_\Gamma : G \onto \Gamma$.  This definition will agree with the definition of $\db_n$ given  for a CW-complex $X$ (as defined in \S3 of \cite{Ha1}) when $G=\pi_1(X)$, $\Gamma=G/{G_r^{(n+1)}}$ and $\phi_\Gamma = \phi_n : G \onto G/{G_r^{(n+1)}}$, the natural projection map. 
For more details see \cite[\S3, \S4 and \S5]{Ha1} and \cite[\S2,\S3,\S5]{Co}.

We recall the definition of a poly-torsion-free-abelian group.  
\begin{definition}\label{ptfa}
A group $\Gamma$ is poly-torsion-free-abelian (PTFA) if it admits
a normal series $\left\{ 1\right\}  =G_{0}\vartriangleleft
G_{1}\vartriangleleft \cdots\vartriangleleft G_{n}=\Gamma$ such
that each of the factors $G_{i+1}/G_{i}$ is torsion-free abelian.
\end{definition}

\begin{remark}
\label{PTFA} Recall that if $A\vartriangleleft G$ is torsion-free-abelian and
$G/A$ is PTFA then $G$ is PTFA. Any PTFA group is torsion-free and
solvable (the converse is not true). Also, any subgroup of a PTFA group
is a PTFA group \cite[Lemma 2.4, p.421]{Pa}.
\end{remark}

Some examples of interesting series associated to a group $G$ are the rational lower central series of $G$
(see Stallings \cite{Sta}), the rational lower central series of the rational commutator subgroup of $G$,
the rational derived series $G_r^{(n)}$ of $G$ (defined below), and the torsion-free derived series $G_H^{(n)}$ of $G$ (see \cite{CH}).
In this paper, our examples and applications will use the rational derived series of a group (defined below).
 We point out that the torsion-free derived series is very interesting since it gives new concordance invariants of
 links in $S^3$ (see \cite{CH} or \cite{Ha2}).  For any of the subgroups $N$ in the above mentioned series, $G/N$ is a PTFA group.
   In particular, for each $n \geq 0$, $\left.G\right/G^{(n+1)}_r$ is PTFA by Corollary 3.6 of \cite{Ha1}.  We recall the definition of $G_r^{(n)}$.
\begin{definition}\label{rds}Let $G$ be a group and $G_{r}^{\left(  0\right)  }=G$.  For $n\geq1$ define
\[  G_{r}^{\left( n\right)  } =\left\{  g\in G_{r}^{\left(  n-1\right) }\mid
g^{k}\in\left[  G_{r}^{\left(  n-1\right)  },G_{r}^{\left(
n-1\right) }\right]  \text{ for some } k \in \Z-\{  0\}
\right\}
\]
to be the \textbf{$n^{th}$ term of the rational derived
series} of $G$.
\end{definition}


R. Strebel showed that if $G$ is the fundamental group of a
(classical) knot exterior then the quotients of successive terms
of the derived series are torsion-free abelian \cite{Str}. Hence
for knot exteriors we have $G_{r}^{\left(  i\right)  }=G^{\left(
i\right)  }$.  This is also well known to be true for free groups.
Since any non-compact surface has free fundamental group, this
also holds for all orientable surface groups.

We make some remarks about PTFA groups. 
Recall that if $\Gamma$ is PTFA then $\Z\Gamma$ is an Ore domain and hence $\Z\Gamma$ embeds in it \emph{right ring of quotients} 
$\mathcal{K}_\Gamma : = \Z\Gamma (\Z\Gamma-\{0\})^{-1}$ which is a skew field.  
More generally, if 
$S\subseteq R$ is a right divisor set of a ring $R$ then the
\emph{right quotient ring} $R S^{-1}$ exists (\cite[p.146]{Pa} or
\cite[p.52]{Ste}).  By $R S^{-1}$ we mean a ring containing $R$
with the property that
\begin{enumerate}
    \item Every element of $S$ has an inverse in $RS^{-1}$.
    \item Every element of $RS^{-1}$ is of the form $rs^{-1}$ with $r \in R$, $s \in S$.
\end{enumerate}
If R is an Ore domain and S is a right divisor set then $RS^{-1}$ is flat as
a left R-module [Ste, Proposition II.3.5]. In particular, $\mathcal{K}_\Gamma$ is a flat left $\Z\Gamma$-module.
Moreover, every finitely generated 
right module over a skew field is free and such modules have a well defined rank, $\rk_{\mathcal{K}_\Gamma}$, 
which is additive on short exact sequences \cite[p.48]{Coh}.  
Thus, 
if C is a non-negative finite chain complex of finitely generated
free right $\Z\Gamma$-modules then the Euler characteristic $\chi(\mathcal{C}) = \sum_{i=0}^\infty (-1)^i \rk C_i$ is
defined and is equal to $\sum_{i=0}^\infty (-1)^i \rk_{\mathcal{K}_\Gamma} H_i(C;\mathcal{K}_\Gamma)$.  In this paper, we will repeatedly use this fact 
about the Euler characteristic.

Let $\psi:G \onto \Z$ be a surjective homomorphism.  Note that we will always be considering $\Z$ as the multiplicative group $\Z=\left<t\right>$ generated by $t$.
We wish to define $\db_\Gamma(\psi)$ as an non-negative integer.  However, in order to do this, we need some compatibility conditions on $\Gamma$ and $\psi$.
\begin{definition}
Let $G$ be a group, $\phi_{_\Gamma}:G \onto \Gamma$, and $\psi:G \twoheadrightarrow \Z$ where  $\Gamma$ is a PTFA group.
 We say that $(\phi_{_\Gamma},\psi)$ is an \textbf{admissible pair} for $G$ if there exists a surjection
$\alpha_{_{\Gamma,\Z}} : \Gamma \twoheadrightarrow \Z$ such that $\psi= \alpha_{_{\Gamma,\Z}} \circ \phi_{_\Gamma}$.
If $\alpha_{_{\Gamma,\Z}}$ is an isomorphism then we say that $(\phi_{_\Gamma},\psi)$ is \textbf{initial}.\end{definition}

Let $(\phi_{_\Gamma},\psi)$ be an admissible pair for $G$.  We define $\Gamma^{\prime} := \ker(\alpha_{_{\Gamma,\Z}})$.
It is clear that $(\phi_{_\Gamma},\psi)$ is initial if and only if $\Gamma^{\prime}=1$.  Since $\Gamma$ is PTFA 
 by Remark~\ref{PTFA}, $\Gamma^{\prime}$ is PTFA.  Hence $\Gamma^{\prime}$ 
embeds in its right ring of quotients which we call $\K_\Gamma$.    Moreover, $\Z\Gamma^{\prime}-\{0\}$ is known to be a 
right divisor set of $\Z\Gamma$ \cite[p. 609]{Pa} hence we can define the right quotient ring 
$R_\Gamma := \Z\Gamma (\Z\Gamma^{\prime}-\{0\})^{-1}$.  After choosing a splitting $\xi: \Z \rightarrow \Gamma$, we see that any
element of $R_\Gamma$ can be written uniquely as $\sum t^{n_i} k_i$ where $t=\xi(1)$ and $k_i \in \K_\Gamma$.  In this way, one sees 
that $R_{\Gamma}$ is isomorphic to the skew polynomial ring $\K_\Gamma[t^{\pm 1}]$ (see the proof of Proposition 4.5 of \cite{Ha1} for more details).
Moreover,  the embedding $g_\psi : \Z\Gamma^{\prime} \rightarrow \K_{\Gamma}$ 
extends to this isomorphism $R_{\Gamma} \rightarrow \K_\Gamma[t^{\pm 1}]$ (here we are identifying $\K_\Gamma$ and $t^0 \K_\Gamma$).   

The abelian group $(G_\Gamma)_{ab}={\ker \phi_{_\Gamma}}\left/{[\ker \phi_{_\Gamma},\ker \phi_{_\Gamma} ]}\right.$ is a right $\Z \Gamma$-module via conjugation,
\[ [g] \gamma = [\gamma^{-1} g \gamma] \]
for  $\gamma \in \Gamma$ and $g \in \ker \phi_{_\Gamma}$.  Moreover, $(G_\Gamma)_{ab}$ is a $\Z\Gamma^{\prime}$-module via the 
inclusion $\Z\Gamma^{\prime} \hookrightarrow \Z\Gamma$.  Thus, $(G_\Gamma)_{ab}\otimes_{\Z\Gamma} \mathcal{K}_\Gamma$ and 
$(G_\Gamma)_{ab}\otimes_{\Z\Gamma^{\prime}} \K_\Gamma$ are right $\mathcal{K}_\Gamma$ and $\K_\Gamma$-modules respectively.
\begin{definition}Let $G$ be a group and $\phi_{_\Gamma}:G \onto \Gamma$ a coefficient
system with $\Gamma$ a PTFA group .  We define the \textbf{$\Gamma$-rank of G} to be 
\[\R_{_\Gamma}(G) = \rk_{\mathcal{K}_\Gamma} \left(\frac{\ker \phi_{_\Gamma}}{[\ker \phi_{_\Gamma},\ker \phi_{_\Gamma} ]} \otimes_{\Z\Gamma} \mathcal{K}_\Gamma \right).\] 
\end{definition}  

For a general group $G$ and coefficient system $\phi_{_\Gamma}$, this rank may be infinite.  However, if $G$ is finitely generated and $\phi_{_\Gamma}$ is 
non-zero then by Proposition~2.11 of \cite{COT}, $r_{_\Gamma}(G)\leq \beta_1(G)-1$ and hence is finite. In the case that 
$\phi_{_\Gamma}$ is the zero map, $r_{_\Gamma}(G)=\beta_1(G)$.

\begin{definition}Let $G$ be a finitely generated group and $(\phi_{_\Gamma},\psi)$ an admissible pair for $G$.  We define the \textbf{$\Gamma$-degree of $\psi$} to be 
\[\db_{\Gamma}(\psi) = \rk_{\K_\Gamma}\left( \frac{\ker \phi_{_\Gamma}}{[\ker \phi_{_\Gamma},\ker \phi_{_\Gamma} ]} \otimes_{\Z\Gamma^{\prime}} \K_\Gamma \right) \] 
if $\R_{_\Gamma}(G)=0$ and $\db_\Gamma(\psi)=0$ otherwise.\end{definition} 

We remark that $(G_\Gamma)_{ab}\otimes_{\Z\Gamma^{\prime}} \K_\Gamma$ is merely $(G_\Gamma)_{ab}\otimes_{\Z\Gamma} \K_\Gamma[t^{\pm 1}]$ viewed as a $\K_\Gamma$-module. 
Since $G$ is a finitely generated group, $(G_\Gamma)_{ab}\otimes_{\Z\Gamma} \K_\Gamma[t^{\pm 1}]$ is a finitely generated 
$\K_\Gamma[t^{\pm 1}]$-module.  Moreover, since $\K_\Gamma[t^{\pm 1}]$ is a (noncommutative left and right) principal ideal domain, \cite[2.1.1, p.49]{Cohn}, the latter is isomorphic to 
\[ \oplus_{i=1}^l \left.\K_\Gamma[t^{\pm 1}]\right/\left<p_i(t)\right> \oplus \left(\K_\Gamma[t^{\pm 1}]\right)^{\R_{_\Gamma}(G)}
\]\cite[Theorem 16, p.43]{Ja}.
Thus, $(G_\Gamma)_{ab}\otimes_{\Z\Gamma^{\prime}} \K_\Gamma$ is a finitely generated $\K_\Gamma$-module if and only if $\R_{_\Gamma}(G)=0$
.  In particular, if $\R_{_\Gamma}(G)= 0$ then $\db_\Gamma(\psi)$ is the sum of the degrees of the $p_i(t)$. Therefore, $\db_\Gamma(\psi)$ as defined above is always finite. 

Let us consider the case when $\Gamma=\Z^m$.  Let $X$ be a CW-complex with $\pi_1(X)=G$ and $X_{\phi_{\Gamma}}$ be the regular $\Z^m$-cover of $X$ corresponding to $\phi_{_\Gamma}$.  Consider an admissible pair $(\phi_{\Z^m},\psi)$ for $G$.  This is one such that $\psi=\psi^{\prime} \circ \phi_\Gamma$ where $\psi^\prime: \Z^m \twoheadrightarrow \Z$.
In this case, $H_1(X_{\phi_{\Gamma}};\Z)= \left.{\ker \phi_{_\Gamma}}\right/{[\ker \phi_{_\Gamma},\ker \phi_{_\Gamma}]}$ is a module over the Laurent polynomial ring with $m$ variables, $\Z[\Z^m]$.  Moreover, $H_1(X_{\phi_{\Gamma}};\Z)$ can be considered as a module over the Laurent polynomial ring with $m-1$ variables $\Z\Gamma^\prime=\Z[\Z^{m-1}]$.  Note that the $m-1$ variables in $\Z[\Z^{m-1}]$ correspond to a choice of basis elements of $\Gamma^\prime=\ker(\alpha_{\Z^m,\Z}: \Z^m \twoheadrightarrow \Z)$. Therefore, as long as the rank of $H_1(X_{\phi_{\Gamma}};\Z)$ as a $\Z[\Z^m]$-module is 0, $\db_{\Z^m}(\psi)$ is equal to the rank of $H_1(X_{\phi_{\Gamma}};\Z)$ as a $\Z[\Z^{m-1}]$-module.  In particular, when $m=1$, $\db_\Z(\psi)$ is equal the rank of $H_1(X_\psi;\Z)$ as an \emph{abelian group} where $X_\psi$ is the infinite cover corresponding to $\psi$ as long as this rank is finite (otherwise $\db_\Z(\psi)=0$).  When $\Z^m$ is the abelianization of $G$, $\db_{\Z^m}(\psi)=\db_0(\psi)$ (see below for the definition of $\db_0$) is equal to the Alexander norm (see \cite{McM} for the definition of the Alexander norm) of $\psi$ by \cite[Proposition~5.12]{Ha1}.

We now define the \emph{higher-order degrees} and \emph{ranks} associated to a group $G$.  For each $n\geq0$, let 
$\Gamma_n=\left.G\right/G_r^{(n+1)}$ where $G_r^{(n+1)}$ is the $(n+1)^{st}$-term of the rational derived series of $G$
as defined in Definition~\ref{rds}.  We 
define the \textbf{$n^{th}$-order rank of $X$} to be \[\R_n(X)= \R_{_{\Gamma_n}}(X).\]
Next, we remark that if $\psi \in H^1(G;\Z)\cong \operatorname{Hom}(G;\Z)$,  then $\psi(G_r^{(n+1)})=1$.  Hence  
for each primitive $\psi \in H^1(G;\Z)$ the pair $(\phi_{_{\Gamma_n}},\psi)$ is an admissible pair for $G$.  
For primitive $\psi$, we define the \textbf{$n^{th}$-order degree of $\psi$} to be \[\db_n(\psi)= \db_{\Gamma_n}(\psi).\]  
For non-primitive $\psi$, there is a primitive cohomology class $\psi^{\prime} \in H^1(X;\Z)$ such that 
$\psi=m \psi^{\prime}$. Define $\db_n(\psi)=m \db_n(\psi^{\prime}).$

Thus, for each group $G$ and $n\geq 0$ we have defined a function $\db_n:H^1(G;\Z) \rightarrow \Z$ which is ``linear on rays through the origin''.  
We put a partial ordering on these functions by  $\db_i \leq \db_j$ if $\db_i(\psi) \leq \db_j(\psi)$ for 
all $\psi \in H^1(G;\Z)$.  Also, we say that $\db_i=0$ provided $\db_i(\psi)=0$ for all $\psi \in H^1(G;\Z)$.  

Suppose $f:E \onto G$ is a surjective homomorphism and $(\phi_{_\Gamma},\psi)$ is an 
admissible pair for $G$. Then there is an induced admissible pair $(\phi_{_\Gamma} \circ f,\psi \circ f)$ for
 $E$.  In particular, we can speak $\db_\Gamma^YE(\psi \circ f)$.  When we have this situation, 
 unless otherwise noted, we will use this admissible pair induced by $G$.  When there is no confusion, we will suppress 
 the $f$ and just write $(\phi_{_\Gamma},\psi )$ when we mean $(\phi_{_\Gamma} \circ f,\psi \circ f)$ or 
 $\psi$ when we mean $\psi \circ f$.

In this paper, we will often use the notation $\R_{_\Gamma}(X)$ and $\db_\Gamma^X(\psi)$ for $X$ a CW-complex and $\psi$ an element of $H^1(X;\Z)\cong H^1(\pi_1(X);\Z)$.  
By this, we mean $\R_{_\Gamma}(\pi_1(X))$ and $\db_\Gamma^{\pi_1(X)}(\psi)$ for an admissible pair $(\phi_{_\Gamma},\psi )$ for $\pi_1(X)$.  These 
are equivalent to the homological definitions given in \cite{Ha1}.  That is, if $(\phi_{_\Gamma},\psi)$ is an admissible pair for $\pi_1(X)$ then  
$H_1(X;\K_{\Gamma}[t^{\pm1}])$ and $H_1(X;\mathcal{K}_\Gamma)$ are right $\K_\Gamma$ and $\mathcal{K}_\Gamma$-modules
respectively and  since $\mathcal{K}_\Gamma$ and $\K_\Gamma[t^{\pm 1}]$ are flat left $\Z\Gamma$-modules \cite[Proposition II.3.5]{Ste}, we see that
\[
\R_{_\Gamma}(X)=\rk_{\mathcal{K}_{\Gamma}}H_1(X;\mathcal{K}_{\Gamma})
\]
and
\[
\db_{\Gamma}(\psi)=\rk_{\K_{\Gamma}} H_1(X;\K_{\Gamma}[t^{\pm1}])
\] 
if $\R_{_\Gamma}(X)=0$ and $\db_{\Gamma}(\psi)=0$ otherwise. 
  
\section{Main Results}\label{mainresults}

We seek to study the behavior of $\db_n(\psi)$ as $n$ increases.  More generally, we would like to compare $\db_\Gamma$ as we vary the group $\Gamma$.  We show that the $\db_\Gamma$ satisfy a monotonicity 
condition provided the groups satisfy a compatibility condition.  We describe this condition below.

\begin{definition}
Let $G$ be a group, $\phi_{_\Lambda}:G \twoheadrightarrow \Lambda$, $\phi_{_\Gamma}:G \twoheadrightarrow \Gamma$, 
and $\psi:G \twoheadrightarrow \Z$ where $\Lambda$ and $\Gamma$ are PTFA groups.  We say that 
$(\phi_{_\Lambda}, \phi_{_\Gamma},\psi)$ is an \textbf{admissible triple} for $G$ if there exist surjections
$\alpha_{_{\Lambda,\Gamma}} : \Lambda \twoheadrightarrow \Gamma$ and 
$\alpha_{_{\Gamma,\Z}} : \Gamma \twoheadrightarrow \Z$ such that $\phi_{_\Gamma}= \alpha_{_{\Lambda,\Gamma}} \circ
\phi_{_\Lambda}$, $\psi= \alpha_{_{\Gamma,\Z}} \circ \phi_{_\Gamma}$, and $\alpha_{_{\Lambda,\Gamma}}$ is not an 
isomorphism.  If $\alpha_{_{\Gamma,\Z}}$ is an isomorphism then we say that $(\phi_{_\Lambda}, \phi_{_\Gamma},\psi)$ 
is \textbf{initial}.
\end{definition}

Note that if $(\phi_{_\Lambda}, \phi_{_\Gamma},\psi)$ an admissible triple then $(\phi_{_\Lambda}, \psi)$ and 
$(\phi_{_\Gamma}, \psi)$ are both admissible pairs. Hence, in this case, we can define both $\db_{\Lambda}(\psi)$ 
and $\db_{\Gamma}(\psi)$.  We note that $(\phi_{_\Lambda}, \phi_{_\Gamma},\psi)$ is initial if and only if 
$(\phi_{_\Gamma}, \psi)$ is initial.  Moreover, $(\phi_{_\Lambda}, \psi)$ is never initial since $\Lambda \onto \Gamma$ 
is not an isomorphism.   We will show that $\db_{\Lambda}(\psi) \geq \db_{\Gamma}(\psi)$ as long as the triple is not 
initial.  We point out that even if $\alpha_{_{\Lambda,\Gamma}}$ is an isomorphism, we can define both the $\Lambda$- and $\Gamma$-degrees and
in this case $\delta_\Gamma(\psi)=\delta_\Lambda(\psi)$!

We now proceed to state and prove the main theorems. 

%

\begin{theorem}\label{2complex}Let $G$ be a finitely presented group with $\Def(G)\geq 1$ and
$(\phi_{_\Lambda},\phi_{_\Gamma},
\psi)$ be an admissible triple for $G$.  If $(\phi_{_\Lambda},\phi_{_\Gamma},\psi)$ is not initial then
\begin{equation}\db_{\Lambda}(\psi) \geq \db_{\Gamma}(\psi)\end{equation} otherwise
\begin{equation} \db_{\Lambda}(\psi) \geq \db_{\Gamma}(\psi)-1.
\end{equation}
\end{theorem}

Before proving Theorem~\ref{2complex}, we will state a Corollary of the theorem and make some remarks about the deficiency hypothesis in the theorem. First,
let $\Gamma_n$ be the quotient of $G$ by the $(n+1)^{st}$ term of the rational derived series as in Definition~\ref{rds}.
Recall that for any $\psi \in H^1(G;\Z)$, $(\phi_{_{\Gamma_n}}, \psi)$ is an admissible pair.  Moreover,
$(\phi_{_{\Gamma_{n+1}}}, \phi_{_{\Gamma_n}}, \psi)$ is an admissible triple unless $G_r^{(n+1)}=G_r^{(n+2)}$ which is initial if and only if
$\beta_1(G)=1$ and $n=0$.  Hence by Theorem~\ref{2complex} we see that the $\db_n$ are a nondecreasing function
of $n$ (for $n \geq 1$).  This behavior was first established for the fundamental groups of knot complements in $S^3$ by T. Cochran in \cite[Theorem 5.4]{Co}.
Recall that $\db_{n+1} \geq \db_n$ (respectively $\db_{n}=0$) means that $\db_{n+1}(\psi) \geq \db_n(\psi)$ (respectively $\db_{n}(\psi)=0$) for all $\psi \in H^1(G;\Z)$.


\begin{corollary}\label{2complexn}Let $G$ be a finitely presented group with $\Def(G)\geq 1$.  
If $\beta_1(G)\geq 2$  then
\[\db_0 \leq \db_1 \leq \cdots \leq \db_n \leq \cdots.
\]
If $\beta_1(G)=1$
and $\psi$ is a generator of $H^1(G;\Z)$ then $\db_0(\psi) -1 \leq \db_1(\psi) \leq \cdots \leq \db_n(\psi) \leq \cdots$.
\end{corollary}


\begin{proof}Let $\psi$ be a primitive class in $H^1(G;\Z)$.  We can assume that $G^{(n+1)}_r \neq G_r^{(n+2)}$ since if $G^{(n+1)}_r =G_r^{(n+2)}$ then
$\db_{n+1}(\psi)=\db_{n}(\psi)$ (note that in the case 
$\beta_1(G)=1$ and $n=0$, $\db_1(\psi)=\db_0(\psi) \geq \db_0(\psi)-1$ is also satisfied).  Therefore  
$T=(\phi_{_{\Gamma_{n+1}}}, \phi_{_{\Gamma_n}}, \psi)$ is 
an admissible triple.  As mentioned above, $T$ is initial if and only if $\beta_1(G)=1$ and $n=0$.  
Hence if $\beta_1(G)=1$ and $n=0$ then by Theorem~\ref{2complex}, $\db_1(\psi) \geq \db_0(\psi)-1$.  
Otherwise, $\db_{n+1}(\psi)\geq\db_n(\psi)$.

If $\beta_1(G)\geq2$ and $\psi$ is not primitive then
$\psi=m \psi^{\prime}$ for some primitive $\psi^{\prime}$ and $m \geq 2$.  Hence,
$\db_{n+1}(\psi)=m \db_{n+1}(\psi^\prime) \geq m \db_{n}(\psi^\prime) = \db_{n}(\psi)$.
\end{proof}

We now make some remarks about the condition $\Def(G)\geq 1$.  First, if $G$ has deficiency 
at least $2$ then the results of Theorem~\ref{2complex} and Corollary~\ref{2complexn} hold 
simply because all of the degrees are zero.

\begin{remark}\label{remark1}If $G$ is a finitely presented group with  $\Def(G)\geq 2$ and $(\phi_{_\Gamma},\psi)$ is an
admissible pair for $G$ then $\R_{_\Gamma}(G)\geq1$ and hence $\db_\Gamma(\psi)=0$.
\end{remark}
To see this, let $X_G$ be a finite, connected 2-complex with one 0-cell $x_0$, $m$ 1-cells, $r$ 2-cells where $m-r\geq 2$ and $G=\pi_1(X_G,x_0)$.  Then $H_1(X_G,x_0;\mathcal{K}_\Gamma)$ has a presentation with $m$
generators and $r$ relations so $\rk_{\mathcal{K}_\Gamma} H_1(X_G,x_0;\mathcal{K}_\Gamma)\geq2$ and hence 
$\R_{_\Gamma}(G)=\R_{_\Gamma}(X_G) = \rk_{\mathcal{K}_\Gamma} H_1(X_G,x_0;\mathcal{K}_\Gamma) -1 \geq 1$ \cite[\S 4 and \S 5]{Ha1}.
Therefore, $\db_\Gamma(\psi)=0$ for all $\psi \in H^1(G;\Z)$.

\vspace{5pt}
However, if the deficiency of G is not positive, we can create an infinite number of examples where the theorem is false!  
We construct finitely presented groups for which the degrees are ``large'' up to (but not
including) the $n^{th}$ stage but the degree at the $n^{th}$ stage is zero!  For simplicity, we only describe examples
when $\beta_1(G)=1$.  However, the reader should notice that the same type of behavior can be seen for groups with $\beta_1(G)\geq 2$ using the same techniques.

\begin{proposition}\label{remark2}For each $g \geq 1$ and $n \geq 1$ there exist examples of finitely presented groups
$G_{n,g}$ with $\Def(G_{n,g})\leq 0$ and $\beta_1(G_{n,g})=1$ such that $\db_0(\psi)=2g$,
$\db_i(\psi)=2g-1$ for $1 \leq i \leq n-1$ and $\db_{n}(\psi)=0$ whenever $\psi$ is a generator of $H^1(G_{n,g};\Z)$.
\end{proposition}
\begin{proof}
We will construct these examples by adding relations to the fundamental group of a fibered knot complement $G$ that kill 
the generators the $\left.G^{(n+1)}\right/G^{(n+2)} \otimes \K_n$.  Let $G$ be the fundamental group of a fibered knot $K$
 in $S^3$ of genus $g \geq 1$ and $n \geq 1$. Since $K$ is fibered, $G^{(1)}$ is free, so $G_r^{(n+1)}/G_r^{(n+2)}=G^{(n+1)}/G^{(n+2)}$ and
 $\mathcal{A}_n=\left.G^{(n+1)}\right/G^{(n+2)} \otimes_{\Z \Gamma_n^{\prime}} \K_n$ is a finitely generated free
 right $\K_n$-module of rank $2g-1$. Let $a_1,\ldots , a_{2g-1}$ be the generators of $\mathcal{A}_n$.  Since $\K_n$ is
 an Ore domain, we can find $k_j \in \K_n$ such that $a_j k_j \in G^{(n+1)}/G^{(n+2)} \otimes 1$. 
 Pick $\gamma_1, \ldots , \gamma_{2g-1} \in G^{(n+1)}$ such that $[\gamma_j]=a_j k_j$ and let
 $H=G/<\gamma_1, \ldots , \gamma_{2g-1}>$ and $\eta: G \onto H$.
Note that since any knot group has deficiency 1, $H$ has a presentation with $m$ generators and $m+2g-2$ relations. Since $\gamma_1, \ldots , \gamma_{2g-1} \in G^{(n+1)}$, we have
an isomorphism $G/G^{(n+1)}\cong H/H^{(n+1)} \cong H/H_r^{(n+1)}$.  Therefore,  $\db_0^H(\psi)=\db_0^G(\psi)=2g$ and $\db_i^Y(\psi)=\db_i^X(\psi)=2g-1$ for $1 \leq i \leq n-1$.

Since $G^\prime \onto H^\prime$, we have $H^{\prime}/H^{(n+1)}\cong G^{\prime}/G^{(n+1)}$ for $0\leq i\leq n$ hence
$\K_n = \K_n^G \cong \K_n^H$.
Moreover, since
$G^{(n+1)} \onto H^{(n+1)}$, the map
$\left.G^{(n+1)}\right/G^{(n+2)} \otimes \K_n \rightarrow \left.H^{(n+1)}\right/H^{(n+2)} \otimes \K_n$ is surjective.
But the generators of $\mathcal{A}_n$ are sent to zero under this map, so $H^{(n+1)}/H^{(n+2)} \otimes \K_n=0$.  
Finally, $H_r^{(n+1)}=H^{(n+1)}$ so  
$$\frac{H_r^{(n+1)}}{H_r^{(n+2)}} \otimes \K_n \cong \frac{H^{(n+1)}}{H_r^{(n+2)}} \otimes \K_n \cong \left(\left.\frac{H^{(n+1)}}{H^{(n+2)}}\right/\{\Z\text{-torsion}\}\right) \otimes \K_n=0$$ 
(see Lemma 3.5 of \cite{Ha1} for the second isomorphism) hence $\db_n(\psi)=0$.
\end{proof}

We will now prove Theorem~\ref{2complex}.

\begin{proof}[Proof of Theorem~\ref{2complex}]If the deficiency of $G$ is at least 2 then by Remark~\ref{remark1}, all of the degrees are zero hence the conclusions of the theorem are true.  Now we prove the case when $\Def(G)=1$.  We can assume that $\R_{_\Gamma}(G) =0$, otherwise $\db_{\Gamma}(\psi)=0$ 
and hence the statement of the theorem is true since $\db_\Lambda(\psi)$ is always non-negative.  
Since $G$ is finitely presented, there is a finite 2-complex $X$ such that $G=\pi_1(X)$ and $\chi(X)=1-\Def(G)=0$.  
Recall that $X$ is obtained from the presentation of $G$ with deficiency 1 by starting with one 0-cell, attaching a 
1-cell for each generator and a 2-cell for each relation in the presentation of $G$.  Since
$\Gamma \twoheadrightarrow \Z$ and $\phi_{_\Gamma}$ is surjective, $H_i(X;\mathcal{K}_\Gamma)=0$ for 
$i \neq 1,2$ \cite[Proposition 2.9]{COT}.  
Moreover, $\chi(X)=0$ implies that  
$\rk_{\mathcal{K}_{\Gamma}} H_2(X;\mathcal{K}_\Gamma)=
\rk_{\mathcal{K}_{\Gamma}} H_1(X;\mathcal{K}_\Gamma)=\R_{_\Gamma}(G)=0$ since the Euler characteristic can be computed using $\mathcal{K}_\Gamma$-coefficients as mentioned in \S-~\ref{definitions}.
Since $\R_{_\Gamma}(X)=0$, it follows that $\R_{_\Lambda}(X)=0$ \cite{Ha2}. 
Replacing $\Gamma$ by $\Lambda$ in the above argument, it follows that  
$\rk_{\mathcal{K}_{\Lambda}} H_2(X;\mathcal{K}_\Lambda)=0$.

Let $X_\psi$ be the infinite cyclic cover of $X$ corresponding to $\psi$.  There is a coefficient system for
$X_\psi$, $\phi^{\prime}_{_\Gamma}:\pi_1(X_\psi) \twoheadrightarrow \Gamma^{\prime}$, given by restricting
$\phi_{_\Gamma}$ to $\pi_1(X_\psi)$. Moreover, as $\K_\Gamma$-modules $H_1(X;\mathcal{K}_\Gamma)\cong H_1(X_\psi;\K_\Gamma)$ so $H_1(X_\psi;\K_\Gamma)$ is a finitely generated free 
$\K_\Gamma$-module of rank $\db_{\Gamma}(\psi)$ (similarly for $\Lambda$).
Since $\Gamma^{\prime}$ is PTFA (and hence $\Z \Gamma^{\prime}$ is an Ore domain), there exists a wedge of $e$ circles 
$W$ and a map 
$f:W \rightarrow X_\psi$ such that 
\[
f_*:H_1(W;\K_\Gamma) \rightarrow H_1(X_\psi;\K_\Gamma)
\]
is an isomorphism.  Here, the coefficient system on $W$ is given by $\phi^{\prime}_{_\Gamma} \circ f_{\ast}$. By the proof of Lemma 2.1 in \cite{COT}, $\ker \phi_{_\Gamma} \neq \ker \psi$ if and only if $\phi^{\prime}_{_\Gamma} \circ f_{\ast}$ is non-trivial.  
Moreover, since $W$ is a finite connected 2-complex with $H_2(W)=0$, if $\ker \phi_{_\Gamma} \neq \ker \psi$ then 
$H_1(W;\K_\Gamma)\cong \K_\Gamma^{e-1}$ \cite[Lemma 2.12]{COT}; otherwise $H_1(W;\K_\Gamma)\cong \K_\Gamma^e$.

Up to homotopy we can assume that $W$ is a subcomplex of $X_\psi$ by replacing $X_\psi$ with the mapping cylinder of $f$.  
Consider the long exact sequence of the pair $(X_\psi,W)$ with coefficients in $\K_\Gamma$: 
\[
H_2(X_\psi;\K_\Gamma) \rightarrow H_2(X_\psi,W;\K_\Gamma) \rightarrow H_1(W;\K_\Gamma) \rightarrow H_1(X_\psi;\K_\Gamma).
\]
Since $X$ has no 3-cells, there is a cell complex, $C_i(X;\Z\Gamma)$, which has no 3-cells.  Therefore, $TH_2(X;\Z\Gamma)$, the $\Z\Gamma$-torsion
submodule of $H_2(X;\Z\Gamma)$, is zero.  Now, the kernel of the map $H_2(X;\Z\Gamma) \rightarrow H_2(X;\mathcal{K}_\Gamma)$ is $TH_2(X;\Z\Gamma)$.  
Moreover, we have shown that $H_2(X;\mathcal{K}_\Gamma)=0$ hence $H_2(X;\Z\Gamma)=0$. Thus, $H_2(X_\psi;\K_\Gamma)
\cong H_2(X;\K_\Gamma[t^{\pm 1}]) \cong H_2(X;\Z\Gamma)\otimes_{\Z\Gamma}\K[t^{\pm 1}] =0$.
Since the last arrow in the sequence is an isomorphism, $H_2(X_\psi,W;\K_\Gamma)=0$.  
Our goal is to show that $H_2(X_\psi,W;\K_\Lambda)=0$. Then by analyzing the long exact sequence of the pair $(X_\psi,W)$ 
with coefficients in 
$\K_\Lambda$, it will follow that $H_1(W;\K_\Lambda) \rightarrow H_1(X_\psi;\K_\Lambda)$ is a monomorphism. 
 We note that $\ker \phi_{_\Lambda} \neq
\ker \phi_{_\Gamma}$ implies that $\rk_{\K_\Lambda} H_1(W;\K_\Lambda)=e-1$ as above.  Thus, if $\ker \phi_{_\Gamma} \neq \ker \psi$
then (assuming the monomorphism above) $\db_\Lambda(\psi) \geq e-1=\db_\Gamma(\psi)$; otherwise  $\db_\Lambda(\psi) \geq e-1=\db_\Gamma(\psi)-1$.

Consider the relative chain complex of $(X_\psi,W)$ with coefficients in $\Z\Gamma^{\prime}$:
\[
0 \rightarrow C_2(X_\psi,W;\Z\Gamma^{\prime}) \xrightarrow{\partial_2^{\Gamma^{\prime}}} C_1(X_\psi,W;\Z\Gamma^{\prime}) \rightarrow.
\]
Since $W$ has no 2-cells, $X_\psi$ has no 3-cells.   Therefore $H_2(X_\psi,W;\Z\Gamma^{\prime})$ is $\Z\Gamma^\prime$-torsion free, so $H_2(X_\psi,W;\K_\Gamma)=0$ implies that 
$H_2(X_\psi,W;\Z\Gamma^{\prime})=0$ 
and hence $\partial_2^{\Gamma^{\prime}}$ is injective.

Let $A=\ker(\alpha_{_{\Lambda,\Gamma}|_{\Lambda^{\prime}}}:\Lambda^{\prime} \twoheadrightarrow \Gamma^{\prime})$. Since $A$ is a subgroup of a PTFA group, A is PTFA by Remark~\ref{PTFA}. 
If $M$ is any right $\Z\Lambda^{\prime}$-module then $M \otimes_{\Z A} \Z$ has the structure of a right $\Z\Gamma^{\prime}$-module given by 
\[(\sum \sigma \otimes n)\gamma = \sum \sigma \gamma \otimes n
\]
for any $\gamma \in \Gamma^{\prime}$.
Moreover, one can check that    
$C_{\ast}(X_\psi,W;\Z\Lambda^{\prime}) \otimes_{\Z A} \Z$ is isomorphic to $C_{\ast}(X_\psi,W;\Z\Gamma^{\prime})$ as right
$\Z\Gamma^{\prime}$-modules.  
Thus, after making this identification, $\partial_2^{\Lambda^{\prime}}: C_2(X_\psi,W;\Z\Lambda^{\prime}) \rightarrow C_1(X_\psi,W;\Z\Lambda^{\prime})$
is injective by the following result of R. Strebel.
\begin{proposition}[R. Strebel, \cite{Str} p. 305]\label{strebel}Suppose $\Gamma$ is a PTFA group and $R$ is a commutative 
ring.  Any map between projective right $R\Gamma$-modules whose image under the functor $- \otimes_{R\Gamma} $ is 
injective, is itself injective. 
\end{proposition}
\noindent Finally, since $\K_{\Lambda}$ is flat as a $\Z\Lambda^{\prime}$-module, $H_2(X_\psi,W;\K_\Lambda)=0$ as desired.
\end{proof}

Suppose $\Lambda$ and $\Gamma$ are abelian groups and $G$ is the fundamental group of a compact orientable manifold with toroidal (or empty) boundary. In this case, it can easily be shown, using the results in \cite{McM} and \cite{Ha1}, that the inequalities in Theorem~\ref{2complex} (and Theorem~\ref{closed} below) are in fact \emph{equalities} for all $\psi$ which lie in the cone of an open face of the Alexander norm ball.
We show below that even in this case, there are $\psi$ for which the inequality in Theorem~\ref{2complex} is necessary.

\begin{example}\label{borr_rings}
Let $X$ be the exterior of the Borromean rings in $S^3$ and let $G$ be the fundamental group of the $X$.  A Wirtinger presentation of $G$ is  given by $\left<x,y,z \hspace{3pt}|\hspace{3pt} [z,[x,y^{-1}]],[y,[z,x^{-1}]]\right>$ (see \cite[p.10]{F} for a similar presentation). Thus, there is an epimorphism $f: G \twoheadrightarrow \left<y,z\right>$ by sending $x$ to $1$. Let $\psi_{(0,m,n)}:G\twoheadrightarrow\Z$ be the homomorphism defined by $\psi(x)=1$, $\psi(y)=t^m$, $\psi(z)=t^n$ where $\text{gcd}(m,n)=1$.  Since $f$ factors through $\psi_{(0,m,n)}$, the rank of $H_1$ of the infinite cyclic cover of $X$ corresponding to $\psi_{(0,m,n)}$ is non-zero (see, for example, \cite[Proposition~2.2]{HaC}). It follows that $\db_\Z(\psi_{(0,m,n)})=0$. However, one can compute the Alexander polynomial of $X$ (from the presentation of $G$) to be $\Delta_X=(x-1)(y-1)(z-1)$.  Therefore, $\db_{\Z^3}(\psi_{(0,m,n)})=|m|+|n|>0$.
\end{example}

Now we consider the case when $G$ is the fundamental group of a \emph{closed} 3-manifold.  In this case, the deficiency of $G$ is 0 so  Theorem~\ref{2complex} does not suffice to prove a monotonicity result for $G$. The proof that the degrees satisfy a monotonicity relation will
use Theorem~\ref{2complex} for 2-complexes but will also use some additional topology of the 3-manifold.  Before stating the
corresponding theorem for closed 3-manifolds, we introduce an  important lemma which will be used in the proof of Theorem~\ref{closed}.

\begin{lemma}\label{whenzero}Let $K$ be a nullhomologous knot in a 3-manifold $X$, $M_K$ be the 0-surgery on $K$, 
$\psi:\pi_1(M_K) \onto \Z$ which maps the meridian of $K$ to a nonzero element of $\Z$, and $(\phi_{_\Gamma},\psi)$ be an 
admissible pair for $\pi_1(M_K)$.  If $\R_{_\Gamma}(M_K)=0$ and $(\phi_{_\Gamma},\psi)$ is not initial then the longitude 
of $K$ is not 0 in $H_1(X \setminus K;\mathbb{K}_\Gamma [t^{\pm 1}])$.  
\end{lemma}
\begin{proof}
Let $l \subset N(K)$ be the longitude of $K$.  Here, $N(K)$ is an open neighborhood of $K$ in $X$.  Note that 
$M_K=(X \setminus N(K)) \cup e^2 \cup e^3$ where the attaching circle of $e^2$ is $l$.  Since $X \setminus N(K)$ is 
homotopy equivalent to $X \setminus K$ we use the latter. Consider the diagram below.  
\[
\begin{diagram} 
\K_\Gamma[t^{\pm 1}] & & & & \\
\dTo_{\partial_3} &  & & & \\
H_2(X \setminus K \cup e^2;\K_\Gamma[t^{\pm 1}]) & \rTo^{\pi} & \K_\Gamma[t^{\pm 1}] & \rTo^{\partial_2} & H_1(X\setminus K; \K_\Gamma[t^{\pm 1}]) \\
\dTo_{i_\ast} & & & & \\
H_2(M_K;\K_\Gamma[t^{\pm 1}])
\end{diagram}
\]
The horizontal (respectively vertical) sequence is the long exact sequence of the pair $\left(X \setminus K \cup e^2,X \setminus K\right)$
(respectively $\left(M_K,X \setminus K \cup e^2\right)$) and the $\K_\Gamma[t^{\pm 1}]$ term in the sequence is generated 
by the relative class coming from $e^2$ (respectively $e^3$).  We note that the boundary of the class represented by $e^2$ 
is the class represented by the longitude of $K$ in $H_1(X\setminus K; \K_\Gamma[t^{\pm 1}])$.  By analyzing the attaching 
map of $\partial e^3$, we see that $\pi \circ \partial_3$ is the map which sends $1$ to $t^r-1$ where $t^r$ is the image 
of the meridian of $K$ under $\phi$.  Since $r \neq 0$ we see that this map is never surjective since $t^r-1$ is not a unit in $\K_\Gamma[t^{\pm 1}]$. 

Since $\R_{_\Gamma}(M_K)=0$, by Remark 2.8 of \cite{COT} we have 
$H_2(M_K;\K_\Gamma[t^{\pm 1}])\cong H^1(M_K;\K_\Gamma[t^{\pm 1}])\cong 
\Ext^1_{\K_\Gamma[t^{\pm 1}]}(H_0(M_K;\K_\Gamma[t^{\pm 1}]),\K_\Gamma[t^{\pm 1}])$.
By the proof of Proposition 2.9 in \cite{COT}, 
$H_0(M_K;\K_\Gamma[t^{\pm 1}]) = \left.\K_\Gamma[t^{\pm 1}]\right/(\K_\Gamma[t^{\pm 1}]\cdot I)$ where $I$ is the 
augmentation ideal of $\Z \pi_1(M_K)$ acting via $\Z\pi_1(M_K) \rightarrow \Z\Gamma \rightarrow \K_\Gamma[t^{\pm 1}]$. 
Thus, $H_0(M_K;\K_\Gamma[t^{\pm 1}])\neq 0$ if and only if $(\phi_{_\Gamma},\psi)$ is initial.  Thus, if 
$(\phi_{_\Gamma},\psi)$ is not initial, $\partial_3$ is surjective.  Suppose $[l]=0$ in $H_1(X\setminus K; \K_\Gamma[t^{\pm 1}])$, 
then $\pi$ would be surjective, making $\pi \circ \partial_3$ surjective which is a contradiction.  
\end{proof}

Consider the situation when $X=S^3\setminus K$, $G=\pi_1(S^3\setminus K)$, $\psi$ is the abelianization map of $G$, and
$\phi_{_\Gamma}:G \twoheadrightarrow \Gamma= G\left/G^{(2)}\right.$ be the quotient map where $G^{(n)}$ is the $n^{th}$ 
term of the derived series of $G$.  It is known that $\Gamma$ is a PTFA group \cite{Str}.  Let $l$ be the longitude of 
$K$.  Since $l \in G^{(2)}$, $\phi_{_\Gamma}$ extends to a map $\pi_1(M_K) \twoheadrightarrow G\left/G^{(2)}\right.$. 
We note that in this case, the pair $(\phi_{_\Gamma},\psi)$ is initial if and only if the Alexander polynomial is 1. 
The longitude being nonzero in $H_1(S^3\setminus K;\mathbb{K}_\Gamma [t^{\pm 1}])$ implies that $l$ is nonzero in 
$H_1(S^3\setminus K;\Z\Gamma)=\frac{G^{(2)}}{G^{(3)}}$.   Hence, if the Alexander polynomial of $K$ is not 1 then 
$l \not \in G^{(3)}$.  This was first proved by T. Cochran in Proposition 12.5 of \cite{Co}.

We now state our main monotonicity theorem for closed 3-manifolds.

\begin{theorem}\label{closed}Let $G$ be the fundamental group of a closed, orientable, connected 3-manifold and $(\phi_{_\Lambda},\phi_{_\Gamma},
\psi)$ be an admissible triple for $G$.  If $(\phi_{_\Lambda},\phi_{_\Gamma},\psi)$ is not initial then
\begin{equation}\db_{\Lambda}(\psi) \geq \db_{\Gamma}(\psi)\end{equation} otherwise
\begin{equation} \db_{\Lambda}(\psi) \geq \db_{\Gamma}(\psi)-2.
\end{equation}
\end{theorem}

As we saw for finitely presented groups with deficiency 1 (Corollary~\ref{2complexn}), for groups of closed 3-manifolds, when $n \geq 1$, the
$\db_n$ are a nondecreasing function of $n$.

\begin{corollary}\label{closedn}Let $G$ be the fundamental group of a  closed, orientable, connected 3-manifold.  If
$\beta_1(G)\geq 2$  then
\[\db_0 \leq \db_1 \leq \cdots \leq \db_n \leq \cdots.
\]
If $\beta_1(G)=1$ and $\psi$ is
a generator of $H^1(G;\Z)$ then $\db_0(\psi) -2 \leq \db_1(\psi) \leq \cdots \leq \db_n(\psi) \leq \cdots$.
\end{corollary}

\begin{proof}[Proof of Theorem~\ref{closed}]
Let $X$ be a closed, orientable, connected 3-manifold with $G=\pi_1(X)$.  We will need the following lemma which is an extension of a lemma of C. Lescop \cite{L}. 

\begin{lemma}\label{closedX}
Let $X$ be a closed, connected, orientable 3-manifold and $\psi:\pi_1(X) \twoheadrightarrow \Z=\left<t\right>$ be a surjective map.  $X$ can be presented
as surgery on an framed link,  $L=\amalg_{i=1}^{\beta_1(x)} L_i$, in a rational homology sphere $R$ such that
\begin{enumerate} \item the components of $L$ are null-homologous in $R$ \label{nullhom}
\item the surgery coefficients on $L_i$ are all 0 \label{zero}
    \item \label{link} $\operatorname{lk}(L_i,L_j)=0$ for $i \neq j$ and
    \item $\psi(\mu_i)=t^{\delta_{1i}}$ when $\mu_i$ is a meridian of $L_i$ and $\delta_{ij}$ is the Kronecker delta.
 \end{enumerate}
\end{lemma}

\begin{proof}[Proof of Lemma]By Lemma 5.1.1 in \cite{L}, $X$ can be obtained by surgery on a framed link $L$ with 
$\beta_1(X)$ components such that (\ref{zero}), (\ref{link}), and (\ref{nullhom}) are satisfied.  Now we note that 
any automorphism of $H_1(X)/\{\Z\text{-torsion}\}\cong \Z^{\beta_1(X)}$ corresponds to a sequence of handleslides and 
reordering or reorienting of the components of $L$.  Moreover, since $\psi$ is a surjective map to $\Z$, there exists 
an automorphism of $H_1(X)/\{\Z\text{-torsion}\}$ that sends the first basis element to $t$ (a generator of $\Z$) and 
the other basis elements to $t^0=1$.  That is, we can do a sequence of handleslides (along with possible reorienting or 
reordering) to get a new link $L^{\prime}$ for which the meridian of the first component maps to $t$ and the other 
meridians map to $1$.  Since the original surgery coefficients and linking numbers of $L$ were $0$, the same is true 
for $L^{\prime}$.  We also note that the components of $L^{\prime}$ are null-homologous in $R$. 
\end{proof}

By Lemma~\ref{closedX} above, 
$X$ can be presented as surgery on a framed link $L=\amalg_{i=1}^{\beta_1(x)} L_i$, 
in a rational homology sphere $R$ such that the first component, $L_1$, has surgery coefficient $0$,
$\operatorname{lk}(L_1,L_i)=0$ for $i \neq 1$ and    
$\psi(\mu_1)=t$ when $\mu_1$ is a meridian of $L_1$.
Let $l$ be the longitude of $L_1$ and 
$Y^{\prime}$ be the space obtained by performing 0-surgery in $R$ on the components $L_2, \dots , L_k$.  Let $Y=Y^{\prime}-N(L_1)$ where $N(L_1)$
is a open neighborhood of $L_1$ in $Y^{\prime}$.
Finally, 
$X^{\prime}=Y\cup_l D^2$ be the space obtained by adding a 2-disk to $Y$ which identifies $\partial D^2$ 
with $l$.

After picking a basepoint in $Y$ (hence in $X$ and $X^{\prime}$), we note that the inclusion map 
induces an isomorphism
$i_{\ast}:\pi_1(X^\prime) \xrightarrow{\cong} \pi_1(X)$.  Thus any coefficient system $\phi_{_\Gamma}$ for $X$
induces a coefficient system for $X^{\prime}$.  Moreover, if $M$ is any $\Z\Gamma$-module then 
$H_1(X;M) \cong H_1(X^{\prime};M)$.  In particular,  $\R_{_\Gamma}(X)=\R_{_\Gamma}(X^{\prime})$ and 
$\db_\Gamma^{X}(\psi) = \db_\Gamma^{X^{\prime}}(\psi)$ for all $\psi \in H^1(X)\cong H^1(X^{\prime})$.
Since $l$ is null-homologous in $Y$, we can identify $H^1(X^{\prime})$ and $H^1(Y)$.  
We define the coefficient systems and admissible pairs for $\pi_1(Y)$ by pre-composing the coefficient systems and admissible pairs for $\pi_1(X)$ with $\pi_1(Y)\rightarrow \pi_1(X)$ induced by the inclusion $Y\subset X$.

We pick the splitting $s:\Z \rightarrow \Gamma$ which sends $t$ to $\phi_{_\Gamma}(\mu_1)$. 
Now we consider the long exact sequence of the pair $(\xp,Y):$
\begin{equation}\label{les_pair}
\rightarrow H_2(\xp,Y;\K_\Gamma[t^{\pm 1}]) \xrightarrow{\partial_2} H_1(Y;\K_\Gamma[t^{\pm 1}]) \rightarrow  
H_1(\xp;\K_\Gamma[t^{\pm 1}]) \rightarrow 0.
\end{equation}
As a $\K_\Gamma[t^{\pm 1}]$-module $H_2(\xp,Y;\K_\Gamma[t^{\pm 1}]) \cong \K_\Gamma[t^{\pm 1}]$
generated by the relative 2-cell $\alpha$.     
Hence as a $\K_\Gamma$-module, $H_2(\xp,Y;\K_\Gamma[t^{\pm 1}])$ is an infinitely generated free module, 
generated by $\alpha t^k$ for $k \in \Z$.  Since the 2-cell is attached along $l$, we have $\partial \alpha = [l]$.  
We note that $l$ and $\mu_1$ live on $\partial N(L_1)$, hence $[l,\mu_1]=1 \in \pi_1(Y)$.  Thus, $[l] (t-1)=0$ 
in $H_1(Y;\K_\Gamma[t^{\pm 1}])$.  Equivalently, $[l]=[l] t^k$ for all k hence the image of $\partial_2$ as a 
$\K_\Gamma$-module has at most one dimension  and is generated by $[l]$.  

Using the same argument as in the first paragraph of Theorem~\ref{2complex}, we can assume that $\R_{_\Gamma}(X)=0$ and 
$\rk_\Gamma H_2(X;\mathcal{K}_\Gamma)=0$.  Since $[l]$ is $t-1$ torsion, the $\partial_2$ map in the long exact sequence 
of the pair $(X^\prime,Y)$ with coefficients in $\mathcal{K}_\Gamma$ is $0$.  
Since $\R_{_\Gamma}(X^{\prime})=\R_{_\Gamma}(X)=0$, we see that $\R_{_\Gamma}(Y)=0$.
By the Theorem in \cite{Ha2}, $\R_{_\Lambda}(Y)=0$.
Thus, $H_1(Y;\K_\Gamma[t^{\pm 1}])$ and 
$H_1(\xp;\K_\Gamma[t^{\pm 1}])$ are finitely generated right $\K_\Gamma$-modules of dimensions $\db^Y_\Gamma(\psi)$ and 
$\db^{\xp}_\Gamma(\psi)=\db^{X}_\Gamma(\psi)$ respectively.

Since $\R_{_\Gamma}(X)=0$, if $(\phi_{_\Gamma},
\psi)$ is not initial then $[l]\neq 0$ in $H_1(Y;\K_\Gamma[t^{\pm 1}])$ by Lemma~\ref{whenzero}.  Also, we note that if $(\phi_{_\Gamma},
\psi)$ is initial, then $\Gamma=\Z$ so all of the meridians except $\psi_1$ lift to the $\Gamma$-cover.  Moreover, since $l$ is
nullhomologous in $Y$, it bounds a surface $F$ in $Y$.  $F$ will lift to the $\Gamma$-cover which implies that $l=0$ in 
$H_1(Y;\K_\Gamma[t^{\pm 1}])$.  Thus, $(\phi_{_\Gamma},\psi)$ is initial if and only if $[l]=0$ in 
$H_1(Y;\K_\Gamma[t^{\pm 1}])$.  Recall that $(\phi_{_\Lambda},\psi)$ is never initial.

Suppose $(\phi_{_\Lambda},\phi_{_\Gamma},
\psi)$ is not initial.  Then $(\phi_{_\Gamma},
\psi)$ is not initial hence  $\db^Y_\Gamma(\psi)=\db^{\xp}_\Gamma(\psi)+1$ (similarly for $\Gamma$).  Since $Y$ is 
homotopy equivalent to a 2-complex with $\chi(Y)=0$, $\db^Y_\Lambda(\psi) \geq \db^Y_\Gamma(\psi)$ by Theorem~\ref{2complex}. 
Therefore 
\[
\db_\Lambda^{^X}(\psi)=\db_\Lambda^{^{\xp}}(\psi)=\db^{^Y}_\Lambda(\psi) -1 \geq \db^{^Y}_\Gamma(\psi) -1 = \db_\Gamma^{^{\xp}}(\psi) = \db_\Gamma^{^X}(\psi).
\] 
Now suppose $(\phi_{_\Lambda},\phi_{_\Gamma},
\psi)$ is initial.  Then $(\phi_{_\Gamma},
\psi)$ is initial  so $\db^Y_\Gamma(\psi)=\db^{\xp}_\Gamma(\psi)$ but $\db^Y_\Lambda(\psi)=\db^{\xp}_\Lambda(\psi)+1$. 
Since $Y$ is homotopy equivalent to
a 2-complex with $\chi(Y)=0$, $\db^Y_\Lambda(\psi) \geq \db^Y_\Gamma(\psi)-1$ by Theorem~\ref{2complex}. 
Therefore 
\[
\db_\Lambda^{^X}(\psi)=\db_\Lambda^{^{\xp}}(\psi)=\db^{^Y}_\Lambda(\psi) -1 \geq (\db^{^Y}_\Gamma(\psi) -1) -1 = \db_\Gamma^{^{\xp}}(\psi) -2 = \db_\Gamma^{^X}(\psi)-2.
\] 
\end{proof}

We point out that there are other higher-order degrees, $\delta_n(\psi)$, for a CW-complex $X$ defined in terms of the
$\K_n[t^{\pm1}]$-torsion submodule of $H_1(X;\K_n[t^{\pm1}])$ (see \cite{Ha1}). These are equal to $\db_n(\psi)$ when $\R_n(X)= 0$.  It would be very interesting to understand 
the monotonicity behavior of these $\delta_n(\psi)$.  In particular, for $n\geq1$ are the $\delta_n(\psi)$ a nondecreasing 
function of $n$?


\section{Applications}

\subsection{Deficiency of a group and obstructions to a group being the fundamental group of a 3-manifold}
Recall that the higher-order ranks and degrees of a CW-complex $X$ only depend on the fundamental group of $X$.  Hence it 
makes sense to talk about the higher-order ranks and degrees of a finitely presented group.  One consequence of the 
theorems in the previous section is that the higher-order degrees give obstructions to a finitely presented group having 
positive deficiency or being the fundamental group of a 3-manifold. 

\begin{proposition}\label{defG}Let $G$ be a finitely presented group and $(\phi_{_\Lambda},\phi_{_\Gamma},\psi)$ be an 
admissible triple for $G$.
\begin{enumerate}\item Suppose $(\phi_{_\Lambda},\phi_{_\Gamma},\psi)$ is not initial.  
If $\db_\Lambda(\psi) < \db_\Gamma(\psi)$ then $\Def(G) \leq 0$ and $G$ cannot be the fundamental group of 
a compact, orientable 3-manifold (with or without boundary).  
\item Suppose $(\phi_{_\Lambda},\phi_{_\Gamma},\psi)$ is initial. 
If $\db_\Lambda(\psi) < \db_\Gamma(\psi) - 1$ then $\Def(G) \leq 0$ and $G$ cannot be the fundamental group 
of a compact, orientable 3-manifold with at least one boundary component which is not a 2-sphere.  In addition,  if 
$\db_\Lambda(\psi) < \db_\Gamma(\psi) - 2$ then $G$ cannot be the fundamental group of a compact, orientable 3-manifold 
(with or without boundary). \end{enumerate}
\end{proposition}
\begin{proof}
First, suppose that $\Def(G)\geq 1$.  Then, by Theorem~\ref{2complex}, $\db_\Lambda(\psi) \geq \db_\Gamma(\psi)$ when $(\phi_{_\Lambda},\phi_{_\Gamma},\psi)$ is not
initial and $\db_\Lambda(\psi) \geq \db_\Gamma(\psi)-1$ when $(\phi_{_\Lambda},\phi_{_\Gamma},\psi)$ is initial.

Now, suppose that $G$ is the fundamental group of a closed, orientable, 3-manifold $X$.   Then, by Theorem~\ref{closed},
$\db_\Lambda(\psi) \geq \db_\Gamma(\psi)$ when $(\phi_{_\Lambda},\phi_{_\Gamma},\psi)$ is not initial and 
$\db_\Lambda(\psi) \geq \db_\Gamma(\psi)-2$ when $(\phi_{_\Lambda},\phi_{_\Gamma},\psi)$ is initial.
Finally, suppose $G$ is the fundamental group of a connected, orientable 3-manifold with boundary.  If at least 1 boundary component is not a 2-sphere then 
$\Def(G)\geq 1$ in which case the paragraph above applies.  Moreover, if all the boundary components $X$ are 2-spheres then $G$ is the
fundamental group of a closed 3-manifold.  
\end{proof}

We point out that Proposition~\ref{defG} is sometimes very easy to use computationally since the groups $\Lambda$ and $\Gamma$ can be taken to be finitely generated free abelian groups.  Using Proposition~\ref{defG}, one can easily prove the well known fact that $\Z^m$ cannot be the group of a compact 3-manifold when $n\geq4$.  

\begin{example}\label{Zm} Consider the initial triple $(\text{id}_{\Z^m},\psi,\psi)$ for $\Z^m$ where $\psi:\Z^m \twoheadrightarrow \Z$ is any surjective map.  Since $\ker(\psi)\cong \Z^{m-1}$, we see that $\db_\Z(\psi)=m-1$. Moreover, since $\ker(\text{id}_{\Z^m})=0$, we see that $\db_{\Z^m}(\psi)=0$.  Therefore, if $m\geq 4$, 
$0=\db_{\Z^m}(\psi) < \db_\Z(\psi) -2 = m-3$.  Thus, by Proposition~\ref{defG}, for $m\geq4$, $\Def(\Z^m)\leq 0$  and $\Z^m$ cannot be the fundamental group of any compact, connected, orientable 3-manifold.  
\end{example}

If we consider the case when the groups $\Gamma$ and $\Lambda$ are quotients of $G$ by the terms of its rational derived series we have the following immediate corollary to Proposition~\ref{defG}.

\begin{corollary}\label{specific_ob}Let $G$ be a finitely presented group.
\begin{enumerate}\item Suppose $\beta_1(G)\geq 2$.  If there exists a $\psi \in H^1(G;\Z)$, and $m,n\in \Z$ such that $n> m\geq 0$ and $\db_{n}(\psi)<\db_m(\psi)$ then $\Def(G)\leq 0$ and $G$ cannot be the fundamental group of a compact, orientable 3-manifold.
\item  Suppose $\beta_1(G)=1$ and $\psi$ in a generator of $H^1(G;\Z)$. 
\begin{enumerate}\item If there exists $m,n\in \Z$ such that $n> m\geq 1$ and $\db_{n}(\psi)<\db_m(\psi)$ then $\Def(G)\leq 0$ and $G$ cannot be the fundamental group of a compact, orientable 3-manifold.  
\item If there exists an $n\in \Z$ such that $n\geq 1$ and $\db_n(\psi)<\db_0(\psi)-1$ then $\Def(G) \leq 0$ and $G$ cannot be the fundamental group of a compact, orientable 3-manifold with at least one boundary component which is not a 2-sphere.  
In addition, if $\db_n < \db_0 - 2$ then $G$ cannot be the fundamental group of a compact, orientable 3-manifold. \end{enumerate}
\end{enumerate}
\end{corollary}

\begin{example} \label{ex}
We saw that the examples $G_{n,g}$ in Proposition~\ref{remark2} satisfy $\db_1(\psi) < \db_0(\psi) - 1 $ when $n=1$
and $g=1$, $\db_1(\psi) < \db_0(\psi) - 2 $ when $n=1$ and $g \geq 2$, and 
$\db_n(\psi) < \db_{n-1}(\psi)$ when $n \geq 2$.
Thus, by Corollary~\ref{specific_ob}, for each $n\geq1$ and $g\geq1$ the groups $G_{n,g}$ in Proposition~\ref{remark2} have
$\Def(G_{n,g})\leq 0$. Moreover, except in the case that $g=1$ and $n=1$, for each $n\geq 1$ and $g\geq1$, the group
$G_{n,g}$ cannot be the fundamental group of a compact, orientable 3-manifold (with or without boundary).  The group
$G_{1,1}$ cannot be the fundamental group of a compact, orientable 3-manifold with at least one boundary
component which is not a 2-sphere.
\end{example}

\subsection{Obstructions to $X \times S^1$ admitting a symplectic structure}

We will show that a consequence of Corollaries~\ref{2complexn} and \ref{closedn} is that the $\db_n(\psi)$ give obstructions to a 4-manifold 
of the form $X \times S^1$ admitting a symplectic structure. 
It is well known that if $X$ is a closed
3-manifold that fibers over $S^1$ then $X \times S^1$ admits a
symplectic structure.  Taubes asks whether the converse is true.

\begin{question}
[Taubes]\label{taubes}Let $X$ be a 3-manifold such that $X\times
S^{1}$ admits a symplectic structure. Does $X$ admit a fibration
over $S^{1}$?
\end{question}

In \cite{Ha1}, we showed that if $X$ is a 3-manifold that fibers over $S^1$ with $\beta_1(X)\geq2$ and $\psi$ representing
the fibration  then $\db_n(\psi)$ is equal to Thurston norm $\|\psi\|_T$ of $\psi$.  This generalized the work of McMullen 
who showed that the Alexander norm  gives a lower bound for the Thurston norm which is an equality when $\psi$ represents 
a fibration. 

\begin{theorem}[\cite{Ha1}]
\label{delbarthm}Let $X$ be a compact, orientable 3-manifold
(possibly with boundary). For all $\psi\in H^{1}(X;\Z)  $ and $n\geq0$
\[
\db_{n}(\psi)  \leq \| \psi \|_T
\]
except for the case when $\beta_{1} (X) =1$, $n=0$,
$X\ncong S^1 \times S^2$, and $X \ncong S^1 \times D^2$. In this
case, $\db_{0} (\psi) \leq \|
\psi\|_{T}+1+\beta_{3}\left( X\right)  $ whenever $\psi$ is
a generator of $H^{1}\left( X;\mathbb{Z}\right) \cong\mathbb{Z}$.
Moreover, equality holds in all cases when $\psi:\pi_{1}(
X) \twoheadrightarrow \mathbb{Z}$ can be represented by a
fibration $X\rightarrow S^{1}$.
\end{theorem}

Using the work of Meng-Taubes and Kronheimer-Mrowka, S. Vidussi
\cite{Vi} has recently given a proof of McMullen's inequality (that the Alexander norm gives a lower bound for the 
Thurston norm of a 3-manifold)
using Seiberg-Witten theory. This generalizes the work of
Kronheimer \cite{K2} who dealt with the case that $X$ is the
0-surgery on a knot. Moreover, Vidussi shows that if $X\times
S^{1}$ admits a symplectic structure (and $\beta_1\left(X\right)
\geq 2$) then the Alexander and Thurston norms of $X$ coincide on
a cone over a face of the Thurston norm ball of $X$, supporting a positive answer to Question~\ref{taubes} asked by Taubes.

\begin{theorem}
[Kronheimer, Vidussi \cite{K2,V,Vi}]\label{vid} Let $X$ be an closed, irreducible
3-manifold such that $X \times S^1$ admits a symplectic structure.
If $\beta_1(X) \geq 2$ there exists a $\psi \in
H^1(X;\mathbb{Z})$ such that $\|\psi\|_A =
\|\psi\|_T$. If $\beta _{1}( X) =1$ then for
any generator $\psi$ of $H^1(X;\mathbb{Z})$,
$\|\psi\|_A = \|\psi\|_T + 2$.
\end{theorem}

In \cite[Theorem~12.5]{Ha1}, we used Vidussi's result and our result that the $\db_n$ give lower bounds for the Thurston norm \cite[Theorem~10.1]{Ha1} to show that the higher-order degrees of a 3-manifold $X$ give algebraic obstructions to a 4-manifold of the form $X\times S^1$ admitting a symplectic structure.  As a result, we were able to show that the closed, irreducible 3-manifolds (with $\beta_1(X)\geq2$) in Theorem 11.1 of \cite{Ha1} have $\db_0 < \db_1 < \cdots <\db_n$  hence cannot admit a symplectic structure.
However, it was still unknown at this time whether Vidussi's Theorem holds if one replaces the Alexander norm with $\db_n$. In \cite[Conjecture~12.7]{Ha1}, we conjectured this to be true.  Since the Alexander norm is equal to $\db_0$, Vidussi's theorem gives us the case when $n=0$.
We will show that Conjecture~12.7 of \cite{Ha1} is true when $n\geq1$.  This is theoretically important since it gives more evidence that the only symplectic 4-manifolds of the form $X\times S^1$ are such that $X$ fibers over $S^1$, supporting a positive answer to the question of Taubes.

\begin{theorem} \label{symplectic}
Let $X$ be a closed, orientable, irreducible
3-manifold such that $X \times S^1$ admits a symplectic structure.
If $\beta_1(X) \geq 2$ there exists a $\psi \in
H^1(X;\mathbb{Z})$ such that \[\db_0(\psi)=\db_1(\psi)= \cdots =
\db_n(\psi) = \cdots =\|\psi\|_T.
\] If $\beta _{1}( X) =1$ then for
any generator $\psi$ of $H^1(X;\mathbb{Z})$,
\[
\db_0(\psi)-2=\db_1(\psi)= 
\cdots 
 \db_n(\psi) = \cdots =\|\psi\|_T.
\]
\end{theorem}

\begin{proof}If $X$ is a closed, orientable, irreducible, 3-manifold with $\beta_1(X)\geq 2$ such that $X \times S^1$ 
admits a symplectic structure then by Theorem~\ref{vid} there exists a $\psi \in
H^1(X;\mathbb{Z})$ such that $\db_0(\psi)=\|\psi\|_A =\|\psi\|_T$. 
By Corollary~\ref{closedn} and Theorem~\ref{delbarthm},  $\db_0(\psi) \leq \db_n(\psi) \leq \|\psi\|_T$ hence 
for all $n\geq0$, $\db_0(\psi)=\db_n(\psi) = \|\psi\|_T$. Similarly, if $\beta_1(X)=1$ then for $\psi$ a generator of 
$H^1(X;\mathbb{Z})$, $\db_0(\psi)-2=\|\psi\|_T$. 
Since $S^1 \times S^2$ is not irreducible, for $n \geq 1$ we have $\db_0(\psi) -2 \leq \db_n(\psi) \leq \|\psi\|_T$ hence 
$\db_0(\psi)-2=\db_n(\psi) = \|\psi\|_T$.
\end{proof}

\subsection{Behavior of the Thurston norm under a continuous map which is surjective on $\pi_1$}

An important problem in 3-manifold topology is determine the behavior of the Thurston norm under  continuous maps $f : X \rightarrow Y$ between 3-manifolds.  It was shown by D. Gabai in \cite{Ga} that if $f$ is a
p-fold covering map then $||f^{\ast}(\psi)||_T=p \, ||\psi||_T$.  Moreover, Gabai showed that if $f$ is a degree $d$ map then $||f^{\ast}(\psi)||_T \geq |d| \, ||\psi||_T$.
These statements were first conjectured by Thurston his original paper on the Thurston norm in \cite[Conjecture 2(b)]{Th}.  We sketch a proof of the latter, since it does not seem to explicitly appear in \cite{Ga}.

\begin{theorem}[Gabai]Let $f: X \rightarrow Y$ be a degree $d$ map between closed, orientable, 3-manifolds. 
 Then for each $\psi \in H^1(Y;\Z)$, $||f^{\ast}(\psi)||_T \geq |d| \, ||\psi||_T$.
\end{theorem}
\begin{proof}
Let $\psi \in H^1(Y;\Z)$ and $F$ be an embedded (possibly disconnected) surface in $X$ such that $[F]$ is
dual to $f^{\ast}(\psi)$ and $\chi_{-}(F)=||f^{\ast}(\psi)||_T$.
Since the following diagram commutes \cite[Theorem 67.2]{Mu}, $[f(F)]=f_{\ast}([F]) = d (\psi \cap \Gamma_Y)$.

\begin{diagram}H^1(X;\Z) & \rTo^{\cap \Gamma_X}_{\cong} & H_2(X;\Z) \\
\uTo_{f^{\ast}} & & \dTo_{f_{\ast}} \\
H^1(Y;\Z) & \rTo^{\cap d\Gamma_Y} & H_2(Y;\Z)
\end{diagram}

By Corollary 6.18 of \cite{Ga}, $||-||_T = x_s(-)$ where $x_s$ is the singular norm.
Hence $|d|\,||\psi ||_T \leq \chi_{-}(f(F)) \leq \chi_{-}(F) = ||f^{\ast}(\psi)||_T$.
\end{proof}
Recall that a degree one map is surjective on $\pi_1$.  Hence one could ask if the existence of a map $f: X \rightarrow Y$ between compact, orientable, 3-manifolds,
that is surjective on $\pi_1$ suffices to guarantee that $||f^{\ast}(\psi)||_T \geq ||\psi||_T$ for all
$\psi \in H^1(Y;\Z)$.  We will give some (algebraic) conditions on $X$ and $Y$ (i.e. that do not depend on the map $f$)
that will guarantee $||f^{\ast}(\psi)||_T \geq ||\psi||_T$.
 
This question was first asked by J. Simon (see Kirby's Problem List \cite[Question 1.12(b)]{Ki}) for knot complements.
Recall that if $K$ is a nontrivial knot in $S^3$ then $H^1(S^3\setminus K;\Z)\cong \Z$ generated by
$\psi$ and $||\psi||_T = 2 g(K)-1$ where $g(K)$ is the genus of $K$.

\begin{jsimon}[J. Simon] If $J$ and $K$ are knots in $S^3$ and $f: S^3\setminus L \rightarrow S^3 \setminus K$ is surjective
on $\pi_1$, is $g(L) \geq g(K)$?
\end{jsimon}

The answer to the above question is known to be yes when $\delta_0(K)=2 g(K)$.  We strengthen this result to the case when $\delta_{n}(K)=2 g(K)-1$ in Corollary~\ref{newgenresult}.  By $\db_n(K)$ we mean $\db_n(\psi)$ for a generator $\psi$ of $H^1(S^3\setminus K;\Z)\cong \Z$. Note that by Theorems~5.4 and~7.1 of \cite{Co}, $$\db_0(K)-1\leq \db_1(K) \leq \cdots \leq \db_n(K) \leq \cdots \leq 2 g(K)-1.$$  Moreover, by Corollary~7.4 of \cite{Co}, there exist knots $K$ for which $\db_0(K)-1 < \db_1(K) < \cdots < \db_n(K)$.  Therefore, the result in Corollary~\ref{newgenresult} is strict generalization of the previously known result.  

Before we state the results concerning the behavior of the Thurston norm under a surjective map on $\pi_1$, we state and prove the following theorem which describes the behavior of $\db_n$ under a surjective map on $\pi_1$.  We only consider the case that $\Def(G)=1$ since if $\Def(G)\geq 2$ then
by Remark~\ref{remark1}, $\R_0(G)\geq 1$.  

\begin{theorem}\label{group_greater}Let $G$ be either $(1)$ a finitely presented group with $\Def(G)=1$ or $(2)$ 
the fundamental group of closed, connected, orientable 3-manifold.  If $P$ is a group with $\beta_1(P)=\beta_1(G)$, $\R_0(G)=0$, and 
$\rho : G \onto P$ is a surjective map
then for each $n \geq 0$ and $\psi \in H^1(P;\Z)$, $$\db_n(\rho^{\ast}(\psi)) \geq \db_n(\psi).$$
\end{theorem}
\begin{proof}  We will first show that the theorem holds for primitive elements of $H^1(Y;\Z)$.  It will then follow for arbitrary elements of
$H^1(Y;\Z)$ since for any
$k\in \Z$, 
$\rho^{\ast}(k \psi)=k \rho^{\ast}(\psi)$, $\db_n(k \rho^{\ast}(\psi))= |k| \db_n(\rho^{\ast}(\psi))$ and  $\db_n(k \psi)= |k| \db_n(\psi)$. Let $\psi$ be a primitive element of $H^1(P;\Z)$, $G_n = G/G_r^{(n+1)}$, and $P_n = P/P_r^{(n+1)}$.
For each $n \geq 0$, we have two coefficient systems for $G$, $\phi^1_n: G \onto G_n$ and $\phi^2_n: G \onto P_n$, defined by $\phi^1_n(g)=[g]$ and $\phi^2_n(g)=[\rho(g)]$.  Note that $\rho$ induces a surjection $\overline{\rho}:G_n \onto P_n$.  Moreover, $\overline{\rho}$ has non-trivial kernel if and only if
$(\phi_n^1,\phi_n^2,\rho^{\ast}(\psi))$ is an admissible triple.

If $\overline{\rho}$ is an isomorphism, then $\db_{G_n}(\rho^{\ast}(\psi))=\db_{P_n}(\rho^{\ast}(\psi))$.  Suppose $\overline{\rho}$ is an not an isomorphism.  
We remark that $(\phi_n^1,\phi_n^2,\rho^{\ast}(\psi))$ is initial if and only if $\beta_1(P)=1$ and $n=0$.  However, since $\rho$ is surjective and $\beta_1(G)=\beta_1(P)$, we have $\overline{\rho}:G_0 \xrightarrow{\cong} P_0$.  Thus $(\phi_n^1,\phi_n^2,\rho^{\ast}(\psi))$ is never inital and hence 
by Theorems \ref{2complex} and \ref{closed}), $\db_{G_n}(\rho^{\ast}(\psi))\geq\db_{P_n}(\rho^{\ast}(\psi))$.

To finish the proof, we will show that $\db_{P_n}(\rho^{\ast}(\psi))=\db_{P_n}(\psi)$.  Since $\rho$ is surjective, we have a surjective map
$$\rho_{\ast}:H_1(G;\Z P_n)=\frac{\ker(\phi^2_n)}{[\ker(\phi^2_n),\ker(\phi^2_n)]} \onto \frac{P_r^{(n+1)}}{[P_r^{(n+1)},P_r^{(n+1)}]}=H_1(P;\Z P_n).$$  
Moreover, since $\K^P_n[t^{\pm 1}]$ 
is a flat (right) $\Z P_n$-module, $\rho_{\ast}:H_1(G;\K^P_n[t^{\pm 1}]) \onto H_1(P;\K^P_n[t^{\pm 1}])$ is surjective. 
The condition $\R_0(G)=0$ implies that both of these modules are torsion \cite{Ha2} hence $\rk_{\K^P_n} H_1(G;\K^P_n[t^{\pm 1}]) \geq \rk_{\K^P_n} H_1(P;\K^P_n[t^{\pm 1}])$ which completes the proof. 
\end{proof}

\begin{corollary}\label{thurston_greater}Suppose there exists an epimorphism $\rho: \pi_1(X) \onto \pi_1(Y)$, where $X$ and $Y$ are compact, connected, 
orientable $3$-manifolds, with toroidal or empty boundaries, such that $\beta_1(X)=\beta_1(Y)$ and $\R_0(X)=0$.  
Let $\psi \in H^1(\pi_1(Y);\Z)$. If any of the following 
conditions is  satisfied 
\begin{description}
    \item[a] $\beta_1(Y) \geq 2$ and $\db_n(\psi)=||\psi||_T$ for some $n \geq 0$
    \item[b] $\beta_1(Y)=1$ and $\db_n(\psi)=||\psi||_T$  for some $n \geq 1$
    \item[c] $\beta_1(Y)=1$, $\beta_3(X) \leq \beta_3(Y)$, $\psi$ is primitive and $\db_0(\psi)=||\psi||_T+1+\beta_3(Y)$
\end{description}
then $$||\rho^{\ast}(\psi)||_T \geq ||\psi||_T.$$
\end{corollary}

\begin{proof}Let $G=\pi_1(X)$ and $P=\pi_1(Y)$.  If $X$ were
$S^1 \times D^2$ or $S^1 \times S^2$ then $\pi_1(X)\cong\Z$ and $\pi_1(Y)\cong \Z$ hence $\db_n(\psi)=0$ for all $n$. Thus, we would be in case \textbf{b} and would have $||\psi||_T=0$ which trivially satisfies the conclusion of the corollary. Therefore, we can assume that $X$ is neither $S^1 \times D^2$ nor $S^1 \times S^2$.
We also remark that since $\R_0(X)=0$, $\Def(\pi_1(X))\leq 1$ by Remark~\ref{remark1}.  Thus, if \textbf{a} or \textbf{b} is satisfied then by Theorem~10.1 of \cite{Ha1},
$\db_n(\rho^\ast(\psi))\leq ||\rho^\ast(\psi)||_T$.  Hence by Theorem~\ref{group_greater} we have 
$||\rho^\ast(\psi)||_T \geq \db_n(\rho^\ast(\psi))\geq \db_n(\psi)=||\psi||_T$.    If \textbf{c} is
satisfied then by Theorem~10.1 of \cite{Ha1} we have $\db_0(\rho^\ast(\psi))\leq ||\rho^\ast(\psi)||_T+1+\beta_3(X)$. Therefore,
$||\rho^\ast(\psi)||_T \geq \db_0(\rho^\ast(\psi))-1-\beta_3(X) \geq \db_0(\psi)-1-\beta_3(Y)=||\psi||_T$.
\end{proof}

We will now discuss the case when $G$ is the fundamental group of a knot complement.

\begin{corollary}\label{knot_greater}If $J$ and $K$ are knots in $S^3$ such that there
exists a surjective homomorphism $\rho: \pi_1(S^3\setminus L) \onto \pi_1(S^3\setminus K)$ then for
each $n\geq 0$, $\db_n(L) \geq \db_n(K)$.
\end{corollary}
\begin{proof} Let $G=\pi_1(S^3\setminus L)$, $P=\pi_1(S^3\setminus K)$, $\psi_P: P \onto P/P^{(1)} \cong \Z$ be the 
abelianization map, and $\psi_G=\psi_P \circ \rho$.  Since $\rho$ is surjective and $\beta_1(S^3-L)=1$, $\psi_G$ is a generator of $H^1(S^3-L;\Z)$.
By \cite[Proposition~2.11 ]{COT}, $\R_0(G)=0$ hence by Theorem~\ref{group_greater}, 
$\db_n(L)=\db_n(\psi_G) \geq \db_n(\psi_P)=\db_n(K)$.
\end{proof}

\begin{corollary}\label{newgenresult}Suppose $J$ and $K$ are knots in $S^3$ such that there 
exists a surjective homomorphism $\rho: \pi_1(S^3\setminus L) \twoheadrightarrow \pi_1(S^3\setminus K)$.  If $\db_0(K)=2 g(K)$ or $\db_n(K)=2 g(K) -1$ for
some $n \geq 1$ then $g(L)\geq g(K)$. 
\end{corollary}
This corollary follows immediately from Corollary~\ref{thurston_greater}.  Instead of omitting any proof,
we will supply a proof which is a simplified version of the proof of Corollary~\ref{thurston_greater}.  
\begin{proof}We can assume that $L$ is not the unknot since $\phi$ is surjective.  If $n\geq 1$ we have $\db_n(L)\leq 2 g(L)-1$ by 
Theorem~7.1 of \cite{Co} or Theorem~10.1 of \cite{Ha1}.  Hence, by Corollary~\ref{knot_greater},    
$2 g(K)-1=\db_n(K) \leq \db_n(L)\leq 2 g(L)-1$.  In the other case, we have $\db_0(L)\leq 2 g(L)$ so
$2 g(K)=\db_0(K) \leq \db_0(L)\leq 2 g(L)$
\end{proof}

\end{document}